\documentclass[12pt]{amsart}
\usepackage{amssymb}
\usepackage{amsmath}
\usepackage{amsfonts,amsthm,latexsym,mathrsfs}
\usepackage[usenames,dvipsnames,svgnames,table]{xcolor}

\theoremstyle{plain}
\newtheorem{theorem}[subsection]{{\bf Theorem}}
\newtheorem{corollary}[subsection]{{\bf Corollary}}
\newtheorem{proposition}[subsection]{{\bf Proposition}}
\newtheorem{lemma}[subsection]{{\bf Lemma}}

\theoremstyle{remark}
\newtheorem{remark}[subsection]{{\it Remark}}
\newtheorem{example}[subsection]{{\it Example}}
\numberwithin{equation}{subsection}

\def \Class {\mathcal S} 
\DeclareMathOperator{\ad}{\rho}
\def \Property {({\star})}

\DeclareMathOperator{\GL}{GL}
\DeclareMathOperator{\Aut}{Aut}
\def \ZZ {\mathbb Z}

\def \FF {\mathbb F}
\def \Dih {Dih}
\DeclareMathOperator{\ATLAS}{ATLAS}
\DeclareMathOperator{\Alt}{Alt}
\DeclareMathOperator{\Sym}{Sym}

\DeclareMathOperator{\PSL}{PSL}
\DeclareMathOperator{\SL}{SL}
\DeclareMathOperator{\PGL}{PGL}

\DeclareMathOperator{\Sp}{Sp}
\DeclareMathOperator{\PSp}{PSp}
\DeclareMathOperator{\Sz}{Sz}

\newcounter{ithmcount}

\makeatletter
\def\@author#1{\g@addto@macro\elsauthors{\normalsize%
    \def\baselinestretch{1}%
    \upshape\authorsep#1\unskip\textsuperscript{%
      \ifx\@fnmark\@empty\else\unskip\sep\@fnmark\let\sep=,\fi
      \ifx\@corref\@empty\else\unskip\sep\@corref\let\sep=,\fi
      }%
    \def\authorsep{\unskip,\space}%
    \global\let\@fnmark\@empty
    \global\let\@corref\@empty
    \global\let\sep\@empty}%
    \@eadauthor={#1}
}
\makeatother

\title[]
{Groups in which every non-nilpotent subgroup is self-normalizing}
\author[C. Delizia]{Costantino Delizia}
\address{University of Salerno, Italy}
\email{cdelizia@unisa.it}
\author[U. Jezernik]{Urban Jezernik}
\address{University of Ljubljana, Slovenia}
\email{urban.jezernik@fmf.uni-lj.si}
\author[P. Moravec]{Primo\v z Moravec}
\address{University of Ljubljana, Slovenia}
\email{primoz.moravec@fmf.uni-lj.si}
\author[C. Nicotera]{Chiara Nicotera}
\address{University of Salerno, Italy}
\email{cnicoter@unisa.it}

\date{}

\begin{document}
\maketitle

\begin{abstract}
We study the class of groups having the property that every non-nilpotent
subgroup is equal to its normalizer. These groups are either soluble or perfect. We
completely describe the structure of soluble groups and finite perfect groups with the above property.
Furthermore, we give some structural information in the infinite perfect case.
\end{abstract}


\section{Introduction}
\label{s:intro}

\noindent A long standing problem posed by Y. Berkovich \cite[Problem 9]{Ber09}
is to study the finite $p$-groups in which every non-abelian subgroup contains
its centralizer.

In \cite{Del16}, the finite $p$-groups which have maximal class or exponent $p$
and satisfy Berkovich's condition are characterized. Furthermore, the infinite
supersoluble groups with the same condition are completely classified. Although
it seems unlikely to be able to get a full classification of finite $p$-groups
in which every non-abelian subgroup is self-centralizing,  Berkovich's problem
has been the starting point for a series of papers investigating
finite and infinite groups
in which every subgroup belongs to a certain family or is self-centralizing.
For instance, in \cite{Del13} and \cite{Del15} locally finite or infinite supersoluble
groups in which every non-cyclic subgroup is self-centralizing are described.

A more accessible version of  Berkovich's problem has been proposed by P.
Zalesskii, who asked to classify the finite groups in which every non-abelian
subgroup equals its normalizer. This problem has been solved in \cite{Del17}.

In this paper we deal with the wider class $\Class$ of groups in which every non-nilpotent
subgroup equals to its normalizer.
All nilpotent groups (and hence all finite $p$-groups) are in
$\Class$. It is also easy to see that groups in $\Class$ are either soluble or
perfect. Further obvious examples of groups in $\Class$
include the minimal non-nilpotent groups (that is,
non-nilpotent groups in which every proper subgroup is nilpotent) and groups in
which every subgroup is self-normalizing. Finite minimal non-nilpotent groups
are soluble, and their structure is well known (see \cite[9.1.9]{Rob96}).
Infinite minimal non-nilpotent groups have been first studied in \cite{New64}
(see also \cite{Cas08} for more recent results). These groups are either
finitely generated or locally finite $p$-groups (\v Cernikov groups or Heineken-Mohamed groups).
Ol'shanskii and Rips (see \cite{Ols89}) showed that there exist
finitely generated infinite simple groups all of whose proper non-trivial
subgroups are cyclic of the same order (the so-called Tarski monsters). On the
other side, groups whose non-trivial subgroups are self-normalizing are periodic
and simple. Furthermore, in the locally finite case they are trivial or of prime
order. Again, infinite examples are the Tarski $p$-groups.

We give a full classification of soluble groups lying in the class $\Class$. We
prove that a non-periodic soluble group belongs to the class $\Class$ if and
only if it is nilpotent (Theorem~\ref{t:nonperiodic_soluble}). Moreover, a
periodic soluble group which is not locally nilpotent lies in the class $\Class$
if and only if it is a split extension of a nilpotent $p'$-group by a cyclic
$p$-group whose structure is completely described in
Theorem~\ref{t:periodic_soluble_notlocnilp}. In particular, this result
characterizes non-nilpotent soluble finite groups in the class $\Class$.
Furthermore, a locally nilpotent soluble group belongs to the class $\Class$ if
and only if it is either nilpotent or minimal non-nilpotent
(Theorem~\ref{t:soluble_locnilp}). Also, an infinite polycyclic group lies in
the class $\Class$ if and only if it is nilpotent (Proposition~\ref{p:polycyclic}).

In the last part of the paper we prove that a finite perfect group lies in the class $\Class$ if and only if it is either isomorphic to the group $\PSL_{2}(2^{n})$ where $2^n-1$ is a prime
number, or to the group $\SL_{2}(5)$ (Theorem~\ref{t:finite_perfect}). Finally, we give some information on the structure of infinite perfect groups in the class $\Class$.

Our notation is mostly standard (see for instance \cite{Ber09} and
\cite{Rob96}). In particular, given any group $G$, we will denote by $Z(G)$ the
center of $G$, by $Z^\infty(G)$ the hypercenter of $G$, by $\Phi(G)$ the
Frattini subgroup of $G$, by $G'$ the commutator subgroup of $G$, and, for all
integers $i\geq1$, by  $\gamma_{i}(G)$ the $i$-th term of the lower central
series of $G$.

\section{General properties of groups in $\Class$}
\label{s:general}

\noindent It is very easy to prove that the class $\Class$ is subgroup and quotient closed.
Furthermore, non-nilpotent groups in $\Class$ are directly indecomposable.

Recall that a group $G$ is said to be perfect if it equals its commutator
subgroup $G'$. Clearly, if $G \in \Class$ then $G$ is perfect or $G'$ is
nilpotent. Hence the groups in $\Class$ are either perfect or soluble.

Suppose now that a cyclic group $\langle x \rangle$ acts on a group $H$ by
means of an automorphism $x$. If a subgroup $L$ of $H$ is invariant with respect
to $\langle x \rangle$, we will write $L \leq_x H$. Consider the induced map
\[
\ad_x \colon H \to H, \quad \ad_x(h) = [x, h] = h^{-x} h.
\]
Clearly, if $H$ is abelian then $\ad_x$ is a homomorphism. We will describe groups belonging to the class $\Class$ based
on the following property of $\ad_x$:
\begin{equation} \label{eq:property}
\forall K \leq_x H, \;
( \exists n\geq1 \colon \ad_x^n(K) = 1 \; \lor \; \langle \ad_x (K) \rangle = K).
\tag*{$\Property$}
\end{equation}

\begin{lemma} \label{l:adMprime}
Let $x$ act on $H$ by means of an automorphism.
Then for every $K \leq_x H$ we have $\langle \ad_x(K) \rangle K' = \ad_x(K) K'$.
\end{lemma}
\proof
Let $h_1, h_2 \in K$. Then $[x, h_1 h_2] = [x, h_2] [x, h_1] [x, h_1,
h_2]$. It follows that $\ad_x(h_1 h_2) \equiv \ad_x(h_1) \ad_x(h_2) \pmod{K'}$.
\endproof

The following easy observations are used in the sequel.

\begin{lemma} \label{l:property}
Let $x$ act on $H$ by means of an automorphism.
\begin{enumerate}
\item The action of $x$ is fixed point free if and only if $\ad_x$ is injective.
\item If $\ad_x$ is injective and $H$ is abelian, then $\Property$ implies that $\ad_x$ is an isomorphism.
\item If $\ad_x$ is injective (or surjective) and $H$ is finite, then $\ad_x$ satisfies $\Property$.
\end{enumerate}
\end{lemma}
\proof
(i)
Note that $\ad_x$ is injective if and only if whenever $[x,h] = 1$ it follows
that $h = 1$. This is precisely the same as $x$ acting fixed point freely on
$H$.

(ii)
Of course we can assume that $H$ is non-trivial. If $\ad_x$ is injective then there is no positive integer $n$ with the property that $\ad_x^n(K) = 1$,
and so $\Property$ implies that $\langle \ad_x (H) \rangle = H$.  If in addition $H$ is abelian then $\ad_x$ is a homomorphism, and hence $\langle \ad_x(H) \rangle = \ad_x(H)$. Therefore $\ad_x$ is an isomorphism.

(iii)
If $H$ is assumed to be finite, then $\ad_x$ is injective if and only if it is
surjective. In this case $\ad_x$ is bijective, and we have that $\ad_x(K) = K$ for all $K \leq_x H$. Thus $\ad_x$ satisfies $\Property$.
\endproof

\begin{lemma}\label{l:nilpotency_from_ad}
Let $G=\langle x\rangle H$, where $H$ is a nilpotent normal subgroup of $G$ generated
by a set $Y$. Suppose that there exists $n\ge 1$ such that $\ad_x^n(y)=1$ for every $y\in Y$.
Then $G$ is nilpotent.
\end{lemma}
\proof
By a theorem of Hall (see for instance \cite[Theorem 2.27]{Rob72}) it suffices to show that $G/H'$ is nilpotent. Thus we may
assume that $H$ is abelian. Consider a commutator of the form
$$c=[x_0,[x_1,[\cdots, [x_{n-1},x_n]]]],$$ where $x_i\in Y\cup\{ x\}$. Without loss of generality we may assume that $x_n\in Y$.
If there exists $i\in \{ 1,2,\ldots ,n-1\}$ such that $x_i\in Y$, then $c\in H'=1$. Thus we are left with the case $x_1=x_2=\cdots =x_{n-1}=x$, but then $c=\ad_x^n(x_n)=1$
by the assumption. As every commutator of length $n+1$ is a product of conjugates of commutators of the above form, this shows that $G$ is nilpotent.
\endproof

\begin{lemma} \label{l:semidirect_center}
Let $G = \langle x \rangle \ltimes H$ be a non-nilpotent group where $x$ has prime
order $p$ and $H$ is nilpotent. Assume that $\ad_x$ has property
$\Property$ and suppose that there exists a subgroup $1 \neq K \leq_x H$ such
that $\ad_x^n(K) = 1$. Then $Z(G) \neq 1$.
\end{lemma}
\proof
As $\langle x \rangle \ltimes K$ is nilpotent by Lemma~\ref{l:nilpotency_from_ad}, it has a non-trivial center. Thus
there exists an element $1 \neq h \in C_K(x)$. Now consider the group $\langle x
\rangle \ltimes Z(H)$. By property $\Property$, we either have $\ad_x(Z(H)) =
Z(H)$ or there is a positive integer $n$ such that $\ad_x^n(Z(H)) = 1$. In the latter case, we
certainly have an element that belongs to $Z(H)$ and commutes with $x$, so that
$Z(G) \neq 1$. Suppose now that $\ad_x(Z(H)) = Z(H)$ holds. By property
$\Property$ we have $$\langle Z(H), h \rangle = \ad_x(\langle Z(H), h \rangle) =
\ad_x(Z(H)) = Z(H),$$ and hence $h \in Z(H)$. Thus we again have $Z(G) \neq 1$.
\endproof

The following proposition shows how property $\Property$ is tightly related to
the class $\Class$.

\begin{proposition} \label{p:class_implies_property}
Let $G = \langle x \rangle \ltimes H$ be a group in $\Class$ with $x^p$
acting trivially on a nilpotent subgroup $H$ for some prime $p$.  Then  $\ad_x$ has
property $\Property$.
\end{proposition}
\proof
Let $K \leq_x H$, and suppose $\langle \ad_x(K) \rangle \subsetneq K$. Consider
the subgroup $$L = \langle x \rangle \ltimes \langle \ad_x(K) \rangle K'$$ of $G$.
As $K$ is nilpotent, it follows that $\langle \ad_x(K) \rangle K' \subsetneq K$ (see for instance \cite[Lemma 2.22]{Rob72}). Now
observe that $L$ is a proper normal subgroup of $\langle x \rangle \ltimes K$.
This implies that $L$ is nilpotent, and so $\ad_x^n(K) = 1$ for some positive integer $n$.
Therefore $\ad_x$ has property $\Property$.
\endproof

Let $p$ any prime number. An abelian group $A$ is said to be $p$-divisible if $A
= pA$.

\begin{lemma} \label{l:p_divisible}
Let $x$ be an automorphism of order $p$ of an abelian group $A$. If $\ad_x$
is surjective, then $A$ is $p$-divisible.
\end{lemma}
\proof
Consider $A$ as a $\ZZ[\langle x \rangle]$-module. In this sense, the operator
$\ad_x$ corresponds to the element $1-x \in \ZZ[\langle x \rangle]$. We
have $(1-x)^p \equiv 0$ modulo $p\ZZ[\langle x \rangle]$, and so the image of
$(\ad_x)^p$ is a subgroup of $p \ZZ[\langle x \rangle] A = pA$.
As $\ad_x$ is assumed to be surjective, it follows that $A = pA$.
\endproof

\begin{corollary}
\label{c:free_abelian}
Let $G = \langle x \rangle \ltimes A$ where $A$ is free
abelian of finite rank, $x$ acts fixed point freely on $A$ and $x^p$ acts trivially on $A$ for
some prime $p$. Then $G$ does lie in the class $\Class$.
\end{corollary}
\proof
By Lemma~\ref{l:p_divisible}, the map $\ad_x$ is not surjective. Since $x$ acts
fixed point freely on $A$, the group $G$ is not nilpotent, and so $\ad_x^n(A)
\neq 1$ for all integers $n$. Whence $\ad_x$ does not have property $\Property$.  It
follows from Proposition~\ref{p:class_implies_property} that $G$ does not belong
to the  class $\Class$.
\endproof

\begin{lemma}
\label{l:semidirect_subgroups_conjugate}
Let $G = \langle x \rangle \ltimes H$ be a periodic non-nilpotent group with
$x^p = 1$ for some prime $p$ and $H$ a nilpotent $p'$-group. Assume that $\ad_x$
has property $\Property$. Then every non-nilpotent subgroup $L \leq G$ is
conjugate to a subgroup of the form $\langle x \rangle \ltimes K$ for some $K
\leq_x H$.
\end{lemma}
\proof
Since $H$ is nilpotent, $L$ is not contained in $H$. It easily follows that $L$ contains an element of the form $x
h$ for some $h \in H$. As $G$ is non-nilpotent, it follows from property
$\Property$ and Lemma~\ref{l:adMprime} that $$h \in H = \langle \ad_x(H)
\rangle \subseteq \ad_x(H) H'.$$ Hence we can write $h = \ad_x(h_1) h'$ with
$h_1 \in H$ and $h' \in H'$. Thus $xh = x^{h_1} h'$. After possibly replacing
$L$ by $L^{h_1^{-1}}$, we can assume that $x h' \in L$ for some $h' \in H'$. As
$\ad_x$ has property $\Property$, we have that either $\ad_x^n(H') = 1$ for some positive integer $n$, or $\langle
\ad_x(H') \rangle = H'$.

In the first case $\langle x \rangle \ltimes H'$ is nilpotent by Lemma~\ref{l:nilpotency_from_ad}. The subgroup $\langle x, h' \rangle$ is finitely generated,
periodic and nilpotent, therefore it is finite. Let $c$ be its  nilpotency class, and let
$p^k$ be the largest power of $p$ that divides $c!$. Set $m$ to be a
positive solution to the congruence system
$$\left\{\begin{aligned}
m &\equiv0 &&\pmod{\frac{|h'|c!\exp(\gamma_2(\langle x, h' \rangle))}{p^k}} \\
m &\equiv1 &&\pmod{p^k}.
\end{aligned}\right.$$
Thus $m(m-1)$ divisible by $c!$, $m$ is divisible by $\exp(\gamma_2(\langle x, h'
\rangle))$ and by $|h'|$, and $m$ is coprime to $p$. In particular, $\binom{m}{i}$ is divisible
by $\exp(\gamma_2(\langle x, h' \rangle))$ for all $1 \leq i \leq c$. By the Hall-Petrescu
formula (see for instance \cite[Appendix 1]{Ber09}) we get
\[
(xh')^m = x^m h'^m g_2^{\binom{m}{2}} \dots g_c^{\binom{m}{c}}
\]
with $g_i \in \gamma_i(\langle x, h' \rangle)$. By the choice of $m$, it follows
that $(xh')^m = x^m$. This element belongs to $L$, and since $x$ is of $p$-power
order, we conclude that $x \in L$.

Consider now the case when $\langle
\ad_x(H') \rangle = H'$. Thus $\langle x \rangle \ltimes H'$ is non-nilpotent, and we can repeat the argument from above with $H$ replaced by $H'$.
Since $H$ is nilpotent, after finitely many steps we find a
conjugate of $L$ that contains $x$. Replacing $L$ by this conjugate we can thus
write $L = \langle x \rangle  \ltimes K$ for $K = L \cap H \leq_x H$.
\endproof

\begin{proposition}
\label{p:periodic_class_iff_property}
Let $G = \langle x \rangle \ltimes H$ be a periodic group with
$x^p$ acting trivially on a nilpotent $p'$-group $H$ for some prime $p$.  Then $G \in
\Class$ if and only if $\ad_x$ has property $\Property$.
\end{proposition}
\proof
If $G \in \Class$, then $\ad_x$ has property $\Property$ by
Lemma~\ref{p:class_implies_property}.

Conversely, assume now that $\ad_x$ has property $\Property$. To prove that $G$
belongs to the class $\Class$, take any non-nilpotent subgroup $L$ of $G$. By
Lemma~\ref{l:semidirect_subgroups_conjugate}, we can assume that $L$ is of the form $L =
\langle x \rangle \ltimes K$ for some $K \leq_x H$. Let us now show that $L$ is
self-normalizing in $G$. To this end, take an element $x^j c \in N_G(L)$. Then
$x^{x^j c} = x^c = x \ad_x(c)$, and so we must have $\ad_x(c) \in K$.
Furthermore, we also have that $c = x^{-j}(x^j c) \in N_H(K)$. Note that
$\langle \ad_x^{-1}(K) \cap N_H(K) \rangle \leq_x H.$ It follows that
$$N_G(L) = \langle x \rangle \ltimes  \langle \ad_x^{-1}(K) \cap N_H(K) \rangle.$$
Now, for any $h_1, h_2 \in \langle \ad_x^{-1}(K) \cap N_H(K) \rangle$, we have
that
$$\ad_x(h_1 h_2) = \ad_x(h_2) \ad_x(h_1) [x, h_1, h_2] \in
K \cdot K \cdot [K, N_H(K)] \subseteq K.$$
Therefore $\ad_x$ maps $\langle \ad_x^{-1}(K) \cap N_H(K) \rangle$ into $K$.
Since $N_G(L)$ is not nilpotent, it follows from property $\Property$ that
\[
\langle \ad_x^{-1}(K) \cap N_H(K) \rangle =
\langle \ad_x (  \langle \ad_x^{-1}(K) \cap N_H(K) \rangle ) \rangle \subseteq
K.
\]
This implies $N_G(L) \leq \langle x \rangle \ltimes K = L$, as required.
\endproof

\begin{remark}
Assume $G = \langle x \rangle \ltimes H$ with $H$ an abelian finite group and
$x$ acting so that $x^p$ acts trivially on $H$. Then $G \in \Class$ if an only
if $\ad_x$ has property $\Property$, which in this case is equivalent by
Lemma~\ref{l:property} to $x$ acting fixed point freely on $H$.  Confer
\cite[Theorem 2.13]{Del17}.
\end{remark}

\begin{remark}
Let $G$ be a periodic group in $\Class$ with a splitting
$\langle x \rangle \ltimes H$, and assume that $x^p = 1$. It might not be the case
that $x$ acts fixed point freely on $H$ (see Example~\ref{example:Q8}). In such
a situation, we have $C_G(x) \cap H \neq 1$. Therefore $C_H(x)$ is a subgroup of
$H$ with $\ad_x(C_H(x)) = 1$. It follows from Lemma~\ref{l:semidirect_center}
that $Z(G) \neq 1$. Now consider the factor group $\langle
x \rangle \ltimes H/Z^{\infty}(G)$. This group is centerless, and so by the
above argument $x$ must act fixed point freely on $H/Z^{\infty}(G)$.
\end{remark}

\begin{example}
\label{example:Q8}
Let $x$ be the automorphism of the quaternion group $Q_8$ given by
$$i \mapsto j,\qquad j \mapsto -k.$$ This is an automorphism of order $3$. Form the semidirect
product $G = \langle x \rangle \ltimes Q_8 \cong \SL_{2}(3)$. It is readily
verified that all proper subgroups of $G$ are nilpotent, and so $G \in \Class$.
Note that $x$ has a non-trivial fixed point on $Q_8$, namely $(-1)^x = -1$. So we
have $Z^{\infty}(G) = \langle -1 \rangle$ with $G/\langle -1 \rangle \cong
\Alt(4)$, and $x$ acts fixed point freely on $Q_8/\langle -1 \rangle \cong C_2
\times C_2$.
\end{example}

\section{Soluble groups in the class $\Class$} \label{s:soluble}
\noindent A Fitting group is one which equals its Fitting subgroup. Thus a Fitting group
is a product of nilpotent normal subgroups, and therefore it is locally
nilpotent. If $G \in \Class$ and $F$ denote the Fitting subgroup of $G$, then
clearly $G$ is a Fitting group or $F$ is nilpotent.

\begin{lemma}
\label{l:nonperfect}
Let $G\in\Class$ be a soluble group, and $F$ the Fitting subgroup of $G$. Then:
\begin{enumerate}
\item $G'\leq F$;
\item $G$ is a Fitting group or $G/F$ has prime order;
\item if $G/G'$ is finitely generated then $G$ is a Fitting group or $G/G'$ is cyclic of prime-power order.
\end{enumerate}
\end{lemma}
\proof The statement (i) is obvious since $G'$ is a nilpotent normal subgroup of $G$.

Suppose now that $G$ is not a Fitting group. Since $G'\leq F$ by (i), the group
$G/F$ is abelian. Let $N/F$ be any proper subgroup of $G/F$. Then $N$ is a
proper normal subgroup of $G$, so $N$ is nilpotent. Thus $N\leq F$. Therefore
$G/F$ has no proper non-trivial subgroups, so it has prime order. Hence (ii) is
proved.

In order to prove (iii), suppose that $G/G'$ is
finitely generated and $G$ is not a Fitting group. Let $M_1/G'$ and $M_2/G'$ be maximal subgroups of $G/G'$.
Since $M_1$ and $M_2$ are proper normal subgroups of $G$, they are nilpotent and
hence contained in $F$. If $M_1\neq M_2$ it follows that $G=M_1M_2 \leq F$, a
contradiction. This means that the finitely generated group $G/G'$ has an unique
maximal subgroup. Therefore $G/G'$ is cyclic of prime-power order.
\endproof

\begin{lemma}
\label{l:locallycyclic}
Let $G\in\Class$ be a non-nilpotent soluble group. Then the quotient group $G/G'$ is a locally cyclic $p$-group for some prime $p$, and $G'=\gamma_3(G)$.
\end{lemma}
\proof
Let $H$ and $K$ be proper normal subgroups of $G$. Then $H$ and $K$ are nilpotent, hence $HK$ is nilpotent and so $HK\neq G$. Thus the result follows by \cite[Theorem 2.12]{New64}.
\endproof

\begin{theorem} \label{t:nonperiodic_soluble}
A non-periodic soluble group belongs to the class $\Class$ if and only if it is nilpotent.
\end{theorem}
\proof
Clearly, if $G$ is nilpotent then $G\in\Class$.

Conversely, let $G\in\Class$ be a non-periodic soluble group, and assume that $G$ is non-nilpotent. Then by Lemma~\ref{l:locallycyclic} the quotient group $G/G'$ is a $p$-group for some prime $p$.

First assume that $G$ is torsion-free. Since $G'$ is nilpotent and $G/G'$ is periodic, it follows by \cite[6.33]{Rob72} that $G$ is nilpotent, again a contradiction.

We are left with the case when the torsion subgroup $T$ of $G$ is non-trivial. Since $G/T$ is a torsion-free soluble group belonging to the class $\Class$, it is nilpotent by the above. As $G/G'$ is a $p$-group it follows that the quotient group $(G/T)/(G/T)'$ is a $p$-group. Hence $G/T$ is a $p$-group (see for instance \cite[5.2.6]{Rob96}). Therefore $G$ is periodic, our final contradiction.
\endproof

Next two results give a complete description of periodic soluble groups
in $\Class$. In particular, our next theorem characterizes finite soluble non-nilpotent groups in $\Class$.

\begin{theorem} \label{t:periodic_soluble_notlocnilp}
Let $G$ be a periodic soluble group, and assume that $G$ is not locally nilpotent. Then $G \in \Class$ if and
only if $G$ splits as $G = \langle x \rangle \ltimes H$, where $\langle x
\rangle$ is a $p$-group for some prime $p$, $H$ is a nilpotent $p'$-group, $x^p$
acts trivially on $H$ and $\ad_x$ has property $\Property$.
\end{theorem}
\proof
If $G$ splits according to the above statement, then it follows from
Propostion~\ref{p:periodic_class_iff_property} that $G$ belongs to the class
$\Class$.

Assume now that $G \in \Class$. Let $x$ be an element of $G$ that does not
belong to the Fitting subgroup $F$. Then $x^p \in F$ for some prime $p$, by Lemma~\ref{l:nonperfect}.
After possibly replacing $x$ by one of its powers, we can assume that the order of $x$ is
a power $p^m$ of $p$. As $G$ is not a Fitting group, we have that $F$ is nilpotent. Hence $F$ is a product of its
Sylow subgroups, say $F = \prod_q S_q$. Note that there is at least one prime $q\ne p$ involved: otherwise $F$ is a $p$-group,
hence $G$ is a $p$-group, but this yields that $G$ is locally nilpotent since it is locally finite, a contradiction. The Sylow subgroups of $F$ are all
characteristic in $F$, so the conjugation action of $x$ preserves them.
Therefore $x$ acts component-wise on $F$.

Note that $x^p \in S_p$. Consider the subgroup $P=\langle x, S_p
\rangle$. Clearly, $P$ is a $p$-group.  Assume that $P\ne\langle x \rangle$, and choose an element $y \in
P\setminus \langle x \rangle$. Set $$J=\left\langle x,y,y^x,y^{x^2},\ldots,y^{x^{p^m-1}}\right\rangle.$$
Since $J$ is a finite $p$-group, and $\langle x\rangle\ne J$, it follows that
$$\langle x\rangle\subsetneq N_J(\langle x\rangle)\subseteq N_P(\langle x\rangle).$$
Hence there exists an element $z\in N_P(\langle x\rangle)\setminus\langle x\rangle$. Thus  $G$ contains the subgroup $\langle x \rangle \ltimes
\prod_{q \neq p} S_q$ that is not self-normalizing. By assumption, this subgroup
must be nilpotent. Since $x$ acts component-wise on $F$, it follows that $G$
itself should be nilpotent. This is a contradiction, from which it follows that
$P = \langle x \rangle$, and so $S_p = \langle x^p
\rangle$. This immediately implies that $G$ splits as $G = \langle x \rangle
\ltimes H$ with $x$ acting component-wise on  $H =
\prod_{q \neq p} S_q$.

Since $x^p \in F$, it commutes with all the $q$-Sylow subgroups of $F$ for $q
\neq p$. As the $p$-Sylow subgroup $S_p$ is cyclic, it follows that $x^p\in Z(G)$.

Finally, let $K \leq_x H$ with $\ad_x^n(K) \neq 1$ for all integers $n$. Therefore the
group $\langle x \rangle \ltimes K$ is non-nilpotent. Consider the group
$\langle x \rangle \ltimes \langle \ad_x(K) \rangle K'$. It is a normal subgroup
of $\langle x \rangle \ltimes K$, so it must either be equal to $\langle x
\rangle \ltimes K$, or else it is nilpotent. The latter case implies that the group
$\langle x, \ad_x(K) \rangle$ is nilpotent, which gives that $\ad_x^n(K) = 1$
for some $n$, a contradiction. Hence we get that $\langle \ad_x(K) \rangle K' = K$, and
since $K$ is nilpotent, it follows that $\langle \ad_x(K) \rangle = K$. Thus
$\ad_x$ has property $\Property$.
\endproof

\begin{corollary} \label{c:dihedral}
Let $n > 2$. The dihedral group $\Dih(n)$ of order $2 n$ belongs to $\Class$
if and only either $n$ is a power of $2$ or $n$ is odd.
\end{corollary}

\begin{theorem} \label{t:soluble_locnilp}
A locally nilpotent soluble group lies in the class $\Class$ if and
only if it is either nilpotent or minimal non-nilpotent.
\end{theorem}
\proof
Clearly, nilpotent and minimal non-nilpotent groups belong to the class $\Class$.

Let $G\in \Class$ be a  periodic soluble group which is locally nilpotent, and assume that $G$ is non-nilpotent. We will prove that $G$ is minimal non-nilpotent. For the sake of contradiction, assume that there exists a proper non-nilpotent subgroup $H$ of $G$. Let $B$ be the last term of the derived series of $G$ which is not contained in $H$. Then $HB$ has the proper non-nilpotent subgroup $H$. Hence without loss of generality we may assume that $G=HB$. Put $L=B\cap H$. Then $B'\leq L$, so $L$ is normal in $B$. Obviously $L$ is normal in $H$, thus $L$ is normal in $G$. The normal series $L<B<G$ can be refined to a (general) principal series of $G$ (see for instance \cite[12.4.1]{Rob96}). Let $W/V$ be any factor of this principal series with $W\leq B$. As $G$ is locally nilpotent, the principal factor $W/V$ is central (see for instance \cite[12.1.6]{Rob96}). Hence $[W,G]\leq V$. This implies that $W\leq N_G(HV)=HV$. Therefore $$W=W\cap HV=(W\cap H)V\leq LV=V.$$ This means $L=B$, a contradiction, and that proves our result.
\endproof

\begin{corollary}
\label{c:periodic_soluble_p}
A locally nilpotent soluble group lying in the class $\Class$ is nilpotent or a $p$-group for some prime $p$.
\end{corollary}
\proof
Let $G\in\Class$ be a locally nilpotent soluble periodic group, and assume that $G$ is non-nilpotent. Then by Theorem~\ref{t:soluble_locnilp} the group $G$ is minimal non-nilpotent, and the result follows by \cite[Lemma 4.2]{New64}.
\endproof

Recall that a group is called just-infinite if it is infinite, but each of its proper quotients is finite. In particular a just-infinite group has no non-trivial finite normal subgroups. The Baer radical of a group $G$ is the subgroup generated by all the cyclic subnormal subgroups of $G$. It has been proved in \cite[Theorem 2]{Wil71} that if $G$ is a just-infinite group with non-trivial Baer radical $A$ then $A$ is free abelian of finite rank and $C_G(A)=A$.

\begin{proposition}
\label{p:polycyclic}
Every infinite polycyclic group in $\Class$ is nilpotent.
\end{proposition}
\proof Assume by a contradiction that the result is false, and let $G\in\Class$ be an infinite polycyclic group which is non-nilpotent. Then Lemma~\ref{l:nonperfect} (iii) ensures that $G/G'$ is cyclic of prime-power order. It easily follows (see for instance \cite[Corollary 2.11]{New64}) that $G'=\gamma_3(G)$.

Let consider the set $\mathcal F$ of all normal subgroups $N$ of $G$ such that $G/N$ is infinite and non-nilpotent. Thus $\mathcal F$ is not empty, as it contains the trivial subgroup. Since $G$ satisfies the maximal condition on subgroups, there exists a maximal element $M\in\mathcal F$. Hence $G/M$ is a non-nilpotent infinite polycyclic group in $\Class$. Moreover, if $G/N$ is any infinite quotient of $G/M$, then the maximality of $M$ implies that either $N=M$ or $G/N$ is nilpotent. Suppose the latter holds. Then there exists a positive integer $t$ such that $\gamma_t(G)\leq N$. Thus $G'\leq N$, a contradiction since $G/G'$ is finite. Therefore, at expense of replacing $G$ by $G/M$, we may assume that $G$ is just-infinite.

Let $F$ denote the Fitting subgroup of $G$. Thus $F$ coincides with the Baer radical of $G$. Then $F$ is free abelian of finite rank by \cite[Theorem~2]{Wil71}. Furthermore, by Lemma~\ref{l:nonperfect} (ii) we may assume that $G/F$ has prime order $p$. Hence, for all $x\in G\setminus F$ we can write $G=\langle x\rangle F$ with $x\notin F$ and $x^p\in F$. We claim that $F/C_F(x)$ is torsion-free. Indeed, suppose $a\in F\setminus C_F(x)$ with $a^n\in C_F(x)$ for some integer $n$. Thus $(a^n)^x=a^n$. Since $F$ is abelian it follows that $(a^xa^{-1})^n=1$. Hence $n=0$ since $F$ is normal and torsion-free. Therefore $F/C_F(x)$ is torsion-free, as claimed. Thus $G/C_F(x)$ is infinite, which implies that $C_F(x)=1$. Hence $x$ acts fixed point freely on $F$. By Corollary~\ref{c:free_abelian}, the group $G$ does not belong to the class $\Class$, our final contradiction.
\endproof

\section{Perfect groups in the class $\Class$} \label{s:perfect}

\begin{lemma}
\label{l:nas}
Let $G \in \Class$ be a finite perfect group, and let $F$ denote its Fitting subgroup.
Then $G/F$ is a non-abelian simple group.
\end{lemma}
\proof
If there is a proper normal subgroup $F \leq M < G$, then $M$ must be nilpotent
since $G \in \Class$, and so $M = F$. Thus $G/F$ is simple. As $G$ is also
assumed to be perfect, $G/F$ is non-abelian.
\endproof

We first classify the finite simple groups in $\Class$.  This is done with the
help of the following lemma.

\begin{lemma} \label{l:simple_max_subs}
Let $G$ be a finite simple group. Then $G$ belongs to $\Class$ if and only if
all of its maximal subgroups belong to $\Class$.
\end{lemma}
\proof
Assume that all maximal subgroups of a finite simple group $G$ belong to
$\Class$, and let $H$ be a non-nilpotent proper subgroup of $G$. As $G$ is
simple, we have $N_G(H) < G$, and so there is a maximal subgroup $M \leq G$ with
$N_G(H) \leq M$. Since $M$ belongs to $\Class$, it follows that $N_G(H) = N_M(H)
= H$, as required.
\endproof

\begin{lemma} \label{l:psl2}
The group $\PSL_2(q)$ belongs to $\Class$ if and only if $q = 2^n$ with
$q - 1$ a prime, or $q \leq 5$.
\end{lemma}
\proof
Suppose that $\PSL_2(q)$ belongs to $\Class$ with $q > 5$, and assume first
that $q$ is odd. This group contains dihedral subgroups of orders
$(q-1)/2$ and $(q+1)/2$ by \cite{Dic01}. Unless $q = 7$, at least one of these does not belong
to $\Class$ by Corollary~\ref{c:dihedral}. Note that $\PSL_{2}( 7)$ has a subgroup isomorphic to $\Sym(4)$, so it does not belong to $\Class$.
Whence we can assume that $q = 2^n$ for some $n \geq 3$.
Now $\PSL_2(q)$ contains a diagonal torus of order $q-1$ acting
fixed point freely on the unipotent subgroup of order $q$. It follows from
Lemma~\ref{l:nonperfect} that the torus must be simple, and so
$q-1$ is either trivial or a prime, as required.
Finally, it follows from \cite[Theorem 2.17]{Del17} that such groups indeed
belong to $\Class$.
\endproof

\begin{proposition} \label{p:simple}
A finite non-abelian simple group belongs to $\Class$ if and only if it is
isomorphic to
$\PSL_{2}( 2^n)$, where $2^n - 1$ is a prime.
\end{proposition}

\proof
We reduce the situation to the case of Lemma~\ref{l:psl2} by
using Lemma~\ref{l:simple_max_subs}.

\begin{itemize}
	\item {\it Alternating groups}. It may be verified readily that  $\Alt(n)$ belongs to $\Class$ if and only $n=5$, since $\Sym(4)$ is contained in $\Alt(n)$ for every $n\ge 6$.
	\item {\it Linear groups $\PSL_n(q)$}. If $n = 2$, this case is covered by
  Lemma~\ref{l:psl2}. If $n \geq 3$, then there is a block embedding of
  $\SL_2(q)$ into $\PSL_n(q)$. The image of this subgroup is normalized
  by the class of a diagonal matrix of the form ${\rm diag}(\alpha, \beta, \gamma, 1, \dots, 1)$. As long as $\alpha \neq \beta$, this diagonal matrix does not
  belong to the image of the embedding of $\SL_2(q)$, and so $\PSL_n(q)$
  does not belong to $\Class$. The only exceptional case is when $|\FF_q^\times| = 1$, i.e.,  $q = 2$,
  in which case either $n = 3$ or $\PSL_n(2)$ contains $\SL_3(2)$ via a block
  diagonal embedding. Both of these groups quotient onto $\PSL_3(2) \cong \PSL_2(7)$, which does not belong to $\Class$.
	\item {\it Symplectic groups $\PSp_{2n}(q)$}. If $n = 1$, then $\PSp_2(q) \cong \PSL_2(q)$ and this is covered above. Now let $n > 1$. Letting
  $W$ be a maximal isotropic subspace of the $2n$-dimensional vector space
  on which $\Sp_{2n}(q)$ acts, the stabilizer of the decomposition $W \oplus W^{\perp}$ is $\GL_n(q) \rtimes C_2$, and so $\PSp_{2n}(q)$ contains $\PGL_n(q) \rtimes C_2$. Therefore these groups do not belong to $\Class$.
	\item {\it Unitary groups and orthogonal groups}. Their associated root systems
  contain a subsystem of type $A_2$, and so they contain subgroups that are
  isomorphic to either $\SL_3(q)$ or $\PSL_3(q)$. None of these belong to
  $\Class$ by above. See \cite{BW97}.
	\item {\it Exceptional Chevalley groups}. We have an inclusion $$G_2(q)\subset F_4(q)\subset E_6(q)\subset E_7(q)\subset E_8(q),$$ and the list of maximal subgroups of $G_2(q)$ in \cite[p. 127]{Wil09} shows that $G_2(q)$, and hence all of the above groups, does not belong to $\Class$.
	\item {\it Steinberg groups ${}^2E_6(q^2)$ and ${}^3D_4(q^3)$}. By \cite[Theorem 4.3]{Wil09}, the group ${}^3D_4(q^3)$ has a maximal subgroup which is isomorphic to $G_2(q^3)$, hence it is not in $\Class$ by the above. Similarly, $F_4(q^2)$ embeds into
	${}^2E_6(q^2)$ by \cite[p. 173]{Wil09}, hence the latter is not in $\Class$.
  \item {\it Suzuki groups $\Sz(q)$}. By \cite[Theorem 4.1]{Wil09}, these contain Frobenius
  groups $C_{q + \sqrt{2q} + 1} \ltimes C_4$ whose Fitting subgroups are of
  index $4$. Such groups do not belong to $\Class$ by Lemma~\ref{l:nonperfect}.
	\item {\it Ree families}. By \cite[Theorem 4.2]{Wil09}, $2\times \PSL_{3}({2n+1})$ is a maximal subgroup of ${}^2G_2(3^{2n+1})$, and $\Sz(2^{2n+1})\wr 2$ is a maximal subgroup of ${}^2F_4(2^{2n+1})$ by \cite[Theorem 4.5]{Wil09}. For the remaining case, ${}^2F_4(2)'$,
	we use $\ATLAS$ \cite{ATLAS} to conclude that this group contains $\Sym(6)$.
	\item {\it Sporadic groups}. Inspection of $\ATLAS$ reveals that each of 26 sporadic groups has a maximal subgroup which is clearly not in $\Class$.
\end{itemize}

\endproof

\begin{proposition}
Let $G\in\Class$ be a perfect non-simple finite group, and let $F$ denote its Fitting subgroup.
Assume that the group $G/F$ contains two elements $a$ and $b$ of distinct prime orders
with the additional property that
$N_{G/F}(\langle a \rangle) \supsetneq \langle a \rangle$ and $N_{G/F}(\langle b \rangle) \supsetneq \langle b \rangle$.
Then the group
$S_p/\Phi(S_p)$ embeds into the Schur multiplier $M(G/F)$,  for every $p$-Sylow subgroup $S_p$ of $F$.
\end{proposition}

\begin{remark}
It is easy to find such elements $a,b$ for the simple groups $\PSL_{2}(2^n)$ that
appear in Proposition~\ref{p:simple}.
One can take $a$ to be an involution (normalized by the Sylow $2$-subgroup of order $2^n$) and $b$ a diagonal matrix of order $q-1$ (normalized
by the class of the flip $\left( \begin{smallmatrix} 0 & 1 \\ 1 & 0 \end{smallmatrix} \right)$).
\end{remark}

\proof
The group $F$ is nilpotent, so we can write $F = \prod_q S_q$ where $S_q$ is a
$q$-group. Now fix a prime $p$ and consider $G_1 = G / \prod_{q \neq p} S_q$.
The Fitting subgroup of $G_1$ is isomorphic to $S_p$. Further, consider the
group $G_2 = G_1 / (S_p / \Phi(S_p))$. The Fitting subgroup $F_2$ of $G_2$ is
an elementary abelian $p$-group, and $G_2$ belongs to the class $\Class$.

Write $S=G/F$. Then $S$ is simple by Lemma~\ref{l:nas}. The group $G_2$ acts on its subgroup $F_2$ by conjugation. There is thus an
induced homomorphism $G_2 \to \Aut(F_2)$.
This homomorphism factors through $F_2$, so we get a homomorphism $$\psi \colon S
\cong G_2/F_2 \to \Aut(F_2).$$ As $S$ is a simple group, we have that either
$\psi$ is injective or trivial. Let us show that $\psi$ must be trivial.

For the sake of contradiction, assume that $\ker \psi = 1$. Since $F_2$ is a
$p$-group, at least one of the elements $a,b$ from the statement of the
Proposition is of order coprime to $p$. Without loss of generality, assume this
element is $a$. Now consider the group $H = \langle a, F_2 \rangle \leq G_2$. By
our assumption on the element $a$, the group $H$ is not self-normalized in
$G_2$. But it is also not nilpotent. Indeed, the element $a$ acts nontrivially
on $F_2$ because $\psi$ is an embedding of $S$ into $\Aut(F_2)$. The order of
$a$ is coprime to $p$, so $\psi$ restricted to $\langle a \rangle$ is a
completely reducible representation of $\langle a \rangle$ on the $GF(p)$-vector
space $F_2$. This representation splits as a sum of $1$-dimensional
representations, and so $a$ is a diagonalizable element in the image of $\psi$.
Being non-trivial, we can not have that $\psi(a) - I$ is a nilpotent matrix,
and so the group $H$ can not be nilpotent. This leads to a contradiciton
with the fact that $G_2 \in \Class$.

We therefore have that $\psi$ is trivial, and so $S$ acts trivially on $F_2$.
This means that $F_2$ is central in $G_2$. Since $G_2$ is also assumed to be
perfect, the extension $1 \to F_2 \to G_2 \to S \to 1$ is a stem central
extension. It follows that $F_2 \cong S_p / \Phi(S_p)$ embeds into $M(S)$.
\endproof

Our next result, together with Proposition~\ref{p:simple}, gives a complete classification of all finite perfect groups in $\Class$.
\begin{theorem}
\label{t:finite_perfect}
A finite perfect group $G$ belongs to the class $\Class$ if and only if it is either
simple or isomorphic to $\SL_{2}(5) \cong C_2 \rtimes \PSL_{2}(5)$.
\end{theorem}
\proof
We have $M(\PSL_{2}(5)) \cong C_2$ and all the other $\PSL$'s have trivial Schur
multipliers. So in this case, the only possibility for a non-simple perfect
group $G$ in $\Class$ is a group whose Fitting quotient is $\PSL_{2}(5)$. Such a
group must have $F$ a $2$-group with cyclic Frattini quotient, so $F$ itself is
cyclic. But now as $G/F$ acts trivially on the Frattini quotient of $F$, it
follows that the image of the homomorphism $G/F \to \Aut(F)$ is a $p$-group
\cite[Exercise 4.4]{Khu98}. Since $G/F$ is a non-abelian simple group, this
implies that $G/F$ must act trivially even on $F$. Hence $F$ is central in $G$.
This implies that $G$ is a stem central extension of $G/F$ by $F$, so it follows
that $|F| = 2$ and $G \cong \SL_{2}(5)$.
\endproof

Now we deal with infinite perfect groups in the class $\Class$.

\begin{lemma}
\label{l:infinitesimple}
An infinite perfect group lying in the class $\Class$ is simple if and only if its Fitting subgroup is trivial.
\end{lemma}
\proof
Let $G\in\Class$ be an infinite perfect group, and let $F$ denote its Fitting subgroup.

First suppose $G$ is simple. If $G=F$ then $G$ is nilpotent, a contradiction since $G$ is perfect. Therefore $F=1$.

Now suppose $F=1$, and let $N$ be any proper normal subgroup of $G$. Since $G\in\Class$, the subgroup $N$ is nilpotent, so $N\leq F$. Therefore $N=1$, and $G$ is simple.
\endproof

\begin{lemma}
\label{l:nofitting}
An infinite perfect group lying in the class $\Class$ cannot be a Fitting group.
\end{lemma}
\proof
Let $G\in\Class$ be an infinite perfect group, and suppose that $G$ is a Fitting group. The group $G$ cannot be minimal non-nilpotent by (see \cite[Proposition 144]{Cas08} and \cite[Corollary 1.4]{Asa00}), so there exists a proper non-nilpotent subgroup $H$ of $G$. Choose $x\in G\setminus H$. Since $G$ is generated by its nilpotent normal subgroups, there exists a normal subgroup $N$ of $G$ such that $N$ is nilpotent and $x\in N$. Hence $N\not\subseteq H$. Let $B$ be the last term of the derived series of $N$ which is not contained in $H$. Put $K=HB$. Then $K\in \Class$ is locally nilpotent and non-nilpotent. Put $L=B\cap H$. Thus $L$ is normal in $K$, and the normal series $L<B<K$ can be refined to a (general) principal series of $K$. As in the proof of Theorem~\ref{t:soluble_locnilp}, all factors of this principal series which lie between $L$ and $B$ are trivial. This means $L=B$, a contradiction.
\endproof

Note that the above shows that the finiteness hypothesis in Lemma~\ref{l:nas} may be omitted.

\begin{proposition}
\label{p:nas}
Let $G \in \Class$ be a perfect group, and let $F$ denote its Fitting subgroup.
Then $G/F$ is a non-abelian simple group.
\end{proposition}
\proof
By Lemma~\ref{l:nas} we may assume that $G$ is infinite. Moreover, by Lemmas~\ref{l:infinitesimple} and \ref{l:nofitting} we may assume that $F$ is a non-trivial proper subgroup of $G$. Clearly $F$ is infinite and contains all proper normal subgroups of $G$.
\endproof

We leave it as an open problem whether or not there exist infinite perfect groups in $\Class$ which are not simple. Note that, if such a group $G$ is locally graded and finitely generated, then $G/F$ is still locally graded (see for instance \cite{Lon95}), and hence it has to be finite. Therefore, by Proposition~\ref{p:simple}, $G/F$ is isomorphic to
$\PSL_{2}( 2^n)$, where $2^n - 1$ is a prime.

\end{document}